# DIGITAL FRACTIONAL ORDER CONTROLLERS REALIZED BY PIC MICROPROCESSOR: EXPERIMENTAL RESULTS

**Ivo PETRÁŠ**[1], **Štefan GREGA**[1] **and Ľubomír DORČÁK**[2]

[1] Steve's Electronic Services
37 – 31255 Upper MacLure Rd.
Abbotsford, B.C., V2T 5N4, Canada
e-mail: {*petras,gsteve*}*@telus.net*

[2] Department of Informatics and Process Control
BERG Faculty, Technical University of Košice
B. Němcovej 3, 042 00 Košice, Slovak Republic
e-mail: *lubomir.dorcak@tuke.sk*

**Abstract:** This paper deals with the fractional-order controllers and their possible hardware realization based on PIC microprocessor and numerical algorithm coded in PIC Basic. The mathematical description of the digital fractional - order controllers and approximation in the discrete domain are presented. Numerical approximation of fractional order differentiation and integration is realized by continued fraction expansion of Al–Alaoui operator. An example of realization of the particular case of digital fractional-order $PI^\lambda D^\delta$ controller is shown and described. The hardware realization is based on the microprocessor PIC16F876 (Microchip Technology) and some external devices such as DA converter, etc. Some possibilities of the practical application of this type of controller and other approaches of the digital and analogue realization are also discussed in the article.

## 1 Introduction

Fractional calculus is a generalization of integration and differentiation to non-integer order fundamental operator $_aD_t^\alpha$, $\alpha \in R$, where $a$ and $t$ are the limits of the operation. The two definitions generally used for the fractional differintegral are the Grunwald - Letnikov definition and the Riemann - Liouville definition [e.g. Oldham et al. 1974, Podlubny 1999].

The use of fractional calculus for modeling physical systems has been widely considered in the last decades [e.g. Oldham et al. 1974]. We can also find works dealing with the application of this mathematical tool in control theory [e.g. Axtell et al. 1990, Podlubny 1999], but these works have usually theoretical character, whereas the number of works, in which a real object is analyzed and the fractional - order controller is designed and implemented in practice, is very small. The main reason for this fact is the difficulty of controller implementation. This difficulty arises from the mathematical nature of fractional operators, which, defined by convolution and implying a non-limited memory, demand hard requirements of processors memory and velocity capacities.

## 2 Fractional order controller

The fractional-order $PI^\lambda D^\delta$ controller (FOC) can be described by the fractional-order differential equation [Podlubny 1999]:

$$u(t) = Ke(t) + T_i \, D_t^{-\lambda} e(t) + T_d \, D_t^\delta e(t), \qquad \left(_0 D_t^\alpha \equiv D_t^\alpha \right) \quad (1)$$

The discrete approximation of the fractional-order controller can be expressed as

$$C(z) = K + T_i \left(\omega(z^{-1})\right)^{-\lambda} + T_d \left(\omega(z^{-1})\right)^\delta, \quad (2)$$

where $\lambda$ is an integral order, $\delta$ is a derivation order, $K$ is a proportional constant, $T_i$ is an integration constant $T_d$ is a derivation constant and $\omega(z^{-1})$ is a generating function.

## 3 Digital realizations of fractional order controller by PIC 16F876

The key point in digital implementation of a FOC is the discretization of the fractional order operators. It is well known that, for interpolation or evaluation purposes, rational functions are sometimes superior to polynomials, roughly speaking, because of their ability to model functions with zeros and poles [Vinagre at al. 2000].

So, for direct discretizing $D^{\pm r}$, (0 < r < 1), is highly interesting to use the Al – Alaoui operator, which is mixed scheme of Euler and Tustin operators, as a generating function [Chen et al. 2002]. Furthermore, for control applications, the obtained approximate discrete-time rational transfer function should be stable and minimum phase. It can be shown that the proposed alternative discretization method, that is, the continued fraction expansion of the Al – Alaoui, enjoy the above desirable properties. By using this method the discrete transfer function approximating fractional-order operators can be expressed as:

$$D^{\pm r}(z) = \left(\frac{8}{7T}\right)^{\pm r} \text{CFE}\left\{\left(\frac{1 - z^{-1}}{1 + z^{-1}/7}\right)^{\pm r}\right\}_{p,q} = \left(\frac{8}{7T}\right)^{\pm r} \frac{P_p(z^{-1})}{Q_q(z^{-1})} \quad (3)$$

where $T$ is the sample period, $P$ and $Q$ are polynomials of degrees $p$ and $q$, respectively, in the variable $z^{-1}$.

The proposed realization of the fractional controller is based on the microprocessor PIC16F876. The block diagram of the mentioned realization is depicted in Fig.1. Setting of the controller parameters is realized via serial port RS232 by personal computer (PC). **Basic features of processor PIC16F876**: 256 data memory bytes, and 368 bytes of user RAM, an integrated 5-channel 10-bit AD converter. Two timers. Precision timing interfaces are done through two CCP modules and two PWM modules. PIC16F876 has 22 I/O pins. The operating speed is 20 MHz and power supply voltage is 5 V (fairly well filtered).

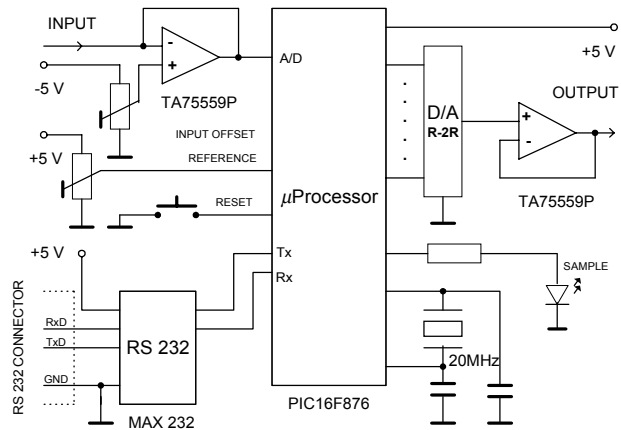

Figure 1. Block diagram of digital FOC based on PIC.

**Basic features of realized digital FOC:** The range of input signal (control error) is 0 - 5 V and can be adjusted for any sensor output. The output signal range (control value) is 0 - 5 V. The reference voltage range for required value setting is also 0 - 5 V. The controller was designed for polynomial degree $p = q = 3$ and for sample period $T \in (0.001 – 120)$ *sec*. This controller was build as a prototype for testing and measurements (see [Petráš at al. 2001].)

## 4 Example: Fractional - order integrator

This work uses the practical case described in [Podlubny at al. 2002]. We realized the fractional half-order integrator as a particular case of fractional order $PI^{\lambda}D^{\delta}$ controller, which has continuous transfer function: $C(s) = 1.4374/s^{0.5}$. The resulting transfer function of half-order integrator was obtained by (3) and for $T = 0.001$ *sec* has the following form:

$$C(z) = 1.4374 \frac{49z^3 - 49z^2 + 7z + 1}{1657z^3 - 2603z^2 + 1048z - 63} \tag{4}$$

The transfer function (4) was rewritten to difference equation and was coded by PIC Basic and then loaded to PIC memory. As the testing signals we used a square impulses (*unit – step* response) and *sinus*. Both signals had amplitude |1V| and frequency 100 Hz. The following figure Fig.2 and Fig.3 presented the measured results (by digital oscilloscope). These results are comparable with results obtained by analog FOC in [Podlubny at al 2002].

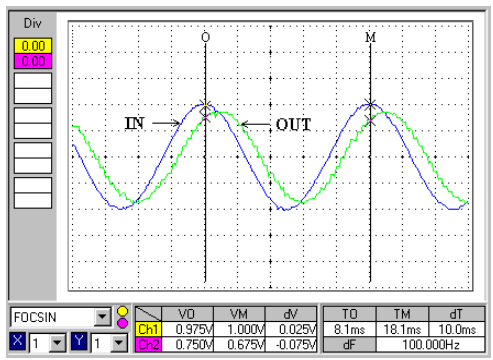
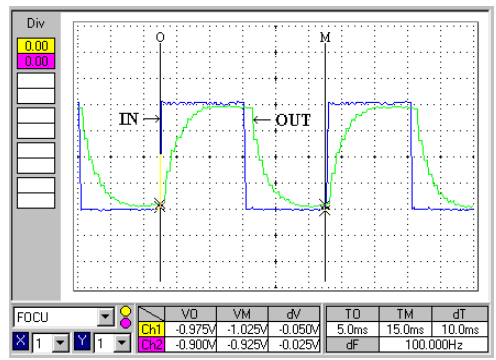

Figure 2. Time response to *sin* input.    Figure 3. Time response to *unit-step* input.

## 5 Conclusion

In this contribution we presented possible realization of digital fractional - order controller based on PIC microprocessor. Approximation of fractional order operator was described via CFE of the Al–Alaoui operator. From the obtained results it can be concluded that for implementing the digital fractional controller is highly interesting to this operator and CFE, because it reduces, without performance degradation, the digital system requirements (e.g. memory and computation time for control law implementation). Such realization of the FOC can be implemented as a single unit device in industrial applications.

In future work we will study the various discretization methods and also other types of PIC processors as for example 18FXX8 series because we had a problem with arithmetic operations in floating point numbers and the speed of mathematical calculation as well. We will try to apply this FOC based on a new type of PIC to control of DC motor speed.

This work was supported by grant VEGA 1/0374/03 from the Slovak Grant Agency.